\documentclass[10pt]{article}

\usepackage{amssymb,amsmath,amsthm,verbatim}
\usepackage{graphicx}
\usepackage{lscape}
\usepackage{indentfirst}
\usepackage{latexsym}
\usepackage{multirow}
\usepackage{multicol}
\usepackage{longtable}
\usepackage{amssymb,amsmath,amsthm,verbatim, graphics}
\usepackage{url}

\newcommand{\bs}{\backslash}
\newcommand{\Z}{\mathbb{Z}}
\newcommand{\Q}{\mathbb{Q}}

\newcommand{\C}{\mathbb{C}}
\newcommand{\h}{\mathfrak{h}}
\newcommand{\R}{\mathbb{R}}

\newcommand{\Pro}{\mathbb{P}}

\newcommand{\GL}{\text{GL}_2}
\newcommand{\M}{\text{M}_2}

\newcommand{\Oh}{\mathcal{O}}

\newcommand{\Xs}{\mathcal{X}^*}
\newcommand{\Gs}{\Gamma^*}

\newcommand{\q}{{\bf q}}

\begin{document}
\title{Minimal Polynomials of Singular Moduli} \author{Eric Errthum} 
\maketitle
\begin{abstract}
Given a properly normalized parametrization of a genus-0 modular curve, the complex multiplication points map to algebraic numbers called singular moduli. In the classical case, the maps can be given analytically. However, in the Shimura curve cases, no such analytical expansion is possible. Fortunately, in both cases there are known algorithms for algebraically computing the rational norms of the singular moduli. We demonstrate a method of using these norm algorithms to algebraically determine the minimal polynomial of the singular moduli below a discriminant threshold. We then use these minimal polynomials to compute the algebraic $abc$-ratios for the singular moduli.
\end{abstract}


\begin{section}{Introduction}
%
%
The classical $j$-function is a modular function that has been studied since the late 1800s by the likes of Gauss, Hermite, Dedekind and Klein. Its values in the upper half plane correspond to isomorphism classes of elliptic curves and  at points corresponding to CM curves, the values of the $j$-function are called singular moduli. A Shimura curve is a generalization of the classical modular curve and again there is a notion of singular moduli and in both cases singular moduli are algebraic numbers. Previously there have been analytical methods for computing singular moduli \cite{Baier} and some limited cases where they can be calculated algebraically \cite{Cox}. In this paper, we present a strictly algebraic method for  determining the defining polynomials for a large class of singular moduli (both classical and Shimura) utilizing pre-existing algorithms that output their rational norms. 

For instance, on the classical modular curve we compute that $$j\left(\frac{1}{2}(1+i\sqrt{39})\right) = -\frac{27}{2} \left(6139474+1702799 \sqrt{13}+147 \sqrt{39 \left(89453213+24809858 \sqrt{13}\right)}\right).$$
Perhaps more importantly, we also compute minimal polynomials of  singular moduli on the Shimura curve of discriminant 6. For example, we find for the singular modulus of discriminant 244, the minimal polynomial is $$x^3 - \frac{2^2 3^7 31\cdot 67\cdot 37223\cdot 235849}{17^6 29^6}x^2 + \frac{2^4 3^{14} 151\cdot 1187\cdot 163327}{17^6 29^6}x -  \frac{2^6 3^{21}19^4}{17^6 29^6}.$$

In Section \ref{Chapter 2} we will give a brief review of modular curves both in the classical and Shimura cases. In Section \ref{Chapter 3} we will review the Gross-Zagier formula for the algebraic norm of the difference of two singular moduli and demonstrate how it can be implemented to algebraically determine the minimal polynomial for a finite collection of singular moduli. In Section \ref{Chapter 4} we adapt this method to the Shimura curve defined by the quaternion algebra of discriminant 6. Finally in Section \ref{Chapter 5} we will review the basics for computing the algebraic $abc$-ratio and present the data obtained from the singular moduli in both cases.

Computations for this paper were performed in Mathematica (v7.0.1.0) and MAGMA (v2.13-14) provided by \cite{Magma}.
\end{section}

\begin{section}{Modular curves and Singular Moduli}\label{Chapter 2}

In this section we provide the basics of Modular Curves and Singular Moduli. For a more complete description of these concepts, see \cite{Ono}, \cite{Shimura}.

The classical modular curve $\Xs_1$ is the one-point compactification of the Riemann surface $\GL(\Z)\bs \h^\pm$ where $\h^\pm = \Pro^1(\C)-\Pro^1(\R)$. As $\Xs_1$ is a genus-$0$ surface, it is isomorphic to $\Pro^1$. Since any map $j: \Xs_1 \to \Pro^1$ must be invariant under the $\GL(\Z)$ action, it suffices to define the function by giving its values at three points. In the classical setting, the points $i=\sqrt{-1}$, $\omega=\frac{1+i\sqrt{3}}{2}$, and $\infty$ are chosen to be sent to $1728$, 0, and $\infty$, respectively. This choice yields what is now commonly referred to as the $j$-function and it has a known Fourier series expansion (where $\q = e^{2\pi i \tau}$)
\begin{eqnarray*}
j(\tau) = \frac{1}{\q}+744+196884\q + \dots \in \frac{1}{\q}\Z[[\q]].
\end{eqnarray*}

Furthermore, the points of $\Xs_1$ correspond to isomorphism classes of elliptic curves. Some of these classes have an extra endomorphism called complex multiplication (CM). Hence, the corresponding points on $\Xs_1$ are called CM points. In the classical case, the CM points of $\Xs_1$ are the imaginary roots of quadratic equations. When $\tau$ is a CM point, $j(\tau)$ is called a singular moduli and is an algebraic integer.

A Shimura curve is a generalization of the modular curve. Let $B$ be the quaternion algebra over $\Q$ with discriminant $D> 1$ and let $\Gamma^* = N_{B^\times}(\Oh) \subset B^\times$ be the normalizer of a maximal order $\Oh \subset B$. Since there is an algebra embedding $B \hookrightarrow \M(\R)$, the discrete group $\Gamma^*$ embeds into $\GL(\R)$ and hence acts on $\h^\pm$. The Shimura curve $\Xs_D$ is then given as 
$$\Xs_D = \Gamma^* \bs \h^\pm.$$ 
When $B$ is a division algebra, $\Xs_D$ is a compact Riemann surface without cusps. 

Points on a Shimura curve can also be identified with certain 2-dimensional abelian varieties and again there is the notion of CM points. 
As before, there will be a generator of the function field, or Hauptmodul, $t: \Xs_D \rightarrow \Pro^1$, and, if properly normalized, the image of a CM point under $t$ will be algebraic over $\Q$. However, since $\Xs_D$ has no cusps, such a map does not have a $\q$-expansion and calculations are more difficult than in the classical case.
\end{section}


\begin{section}{Minimal polynomials of classical singular moduli} \label{Chapter 3}
Let $\tau_1$ and $\tau_2$ be CM points on $\Xs_1$ with relatively prime negative discriminants $-d_1$ and $-d_2$. Then $j(\tau_1)-j(\tau_2)$ is an algebraic integer over $\Q$ of degree $h(-d_1)h(-d_2)$ where $h(-d_i)$ is the class number of the quadratic order of discriminant $-d_i$. Further, let $\omega_i$ be the number of roots of unity in that order. The main result of \cite{GZ} is that the absolute norm $N(\tau_1,\tau_2) :=|j(\tau_1)-j(\tau_2)|_\Q$ is given by
\begin{eqnarray}|j(\tau_1)-j(\tau_2)|_\Q^{\frac{2}{\omega_1\omega2}} =\prod_{\footnotesize{\begin{array}{c}x,n,n' \in \mathbb{Z}\\n,n'>0\\x^2+4nn'=d_1d_2\end{array}}}\!\!\!n^{\epsilon(n')}\label{GZ}\end{eqnarray}
where $\epsilon(p)$ is a completely multiplicative function defined on primes $p$ with $\left(\frac{d_1d_2}{p}\right)\neq 1$ by
$$\epsilon(p) = \left\{\begin{array}{ll}
\left(\frac{-d_1}{p}\right) & \text{if $\gcd(p,-d_1) = 1$,}\\
\left(\frac{-d_2}{p}\right) & \text{if $\gcd(p,-d_2) = 1$.}\end{array}\right.$$

The computational benefit of this method is that very little actually needs to be known about the algebraic integers $j(\tau_i)$ to perform the calculation. 

On the other hand, recall that in general the norm of an algebraic number is given by
$$|\alpha|_\Q = \prod_{\sigma \in \text{Gal}(\Q(\alpha)/\Q)} \sigma(\alpha).$$
Suppose that $r \in \Q$, then $\Q(r-j(\tau)) = \Q(j(\tau))$. Let $G= \text{Gal}(\Q(r - j(\tau))/\Q) = \text{Gal}(\Q(j(\tau))/\Q)$. Then
\begin{align*}
|r-j(\tau)|_\Q &= \pm\prod_{\sigma \in G} \sigma(r-j(\tau))\\
&=\pm\prod_{\sigma \in G} r-\sigma(j(\tau))\\
&=\pm M_{j(\tau)}(r)
\end{align*}
where $M_{\alpha}(x) \in \Q[x]$ denotes the minimal polynomial of the algebraic number $\alpha$.

Since the degree of $M_{j(\tau)}(r)$ is $h(-d_2)$, it suffices to know (up to an issue of sign which we discuss in the example below) the value of the left-hand-side for $h(-d_2)+1$ values of $r$ to interpolate the polynomial. Hence, for a given $j(\tau_0)$ with discriminant $-d_0$, it suffices to find $h(-d_0)+1$ rational $j(\tau_i)$ with $\gcd(d_i,d_0)=1$ and then use (\ref{GZ}) to find the value of $\pm M_{j(\tau_0)}(j(\tau_i))$. 

Unfortunately, the number of rational singular moduli is finite: They occur exactly at the CM points of discriminant $-4$, $-8$, $-3$, $-7$, $-11$, $-19$, $-43$, $-67$, and $-163$. This fact limits the previous method to only being possible for singular moduli of degree 8 or less. The largest achievable case is the singular moduli of discriminant $-5923$.

Also, it should be noted that this algorithm is in no way optimal or the most efficient from a computational viewpoint. Many of these facts can be discovered through much quicker analytic, floating-point calculations and recognition of approximations. The importance of our method lies in it purely algebraic nature. This, in turn, allows it it be  implemented in the Shimura curve case, as we'll see in section 4, where there is no known analytical methods.

\begin{subsection}{Example: $j(\frac{1}{2}(1+i\sqrt{39})$}

Fix $\tau = \frac{1}{2}(1+i\sqrt{39})$. Then $d = 39$ and $h(-d) = 4$. Since $-d$ is relatively prime to the discriminants $-4$, $-8$, $-7$, $-11$, and $-19$, we can use (\ref{GZ}) to compute the following
\begin{eqnarray*}
|j(\tau)-j(1+i)| &=& 3^{12}7^{ 8}19^{4}23^{2}\\
|j(\tau)-j(1+i\sqrt{2})| &=& 7^8 13^2 23^2 29\cdot31^2 37^2 53\\
|j(\tau)-j({\textstyle \frac{1}{2}}(1+i\sqrt{7}))| &=& 3^{12}7^4 13^2 17^3 19^2 31^2\\
|j(\tau)-j({\textstyle \frac{1}{2}}(1+i\sqrt{11}))| &=& 7^8 13^2 17^3 19^2 29^2 101 \cdot107\\
|j(\tau)-j({\textstyle \frac{1}{2}}(1+i\sqrt{19}))| &=& 3^{12} 13^2 19^2 29\cdot31^2 37^2 \cdot53\cdot113\cdot173\cdot179
\end{eqnarray*}
This gives the following 5 points on the curve $y=|M_{j(\tau)}(x) |$:
\begin{align*}
(x_1,y_1) = (12^3,\ &3^{12}7^{ 8}19^{4}23^{2}),\\
 (x_2,y_2) = (20^3,\ &7^8 13^2 23^2 29\cdot31^2 37^2 53),\\
 (x_3,y_3) = (-15^3,\ & 3^{12}7^4 13^2 17^3 19^2 31^2),\\
 (x_4,y_4) = (-32^3,\ &7^8 13^2 17^3 19^2 29^2 101 \cdot107),\\
 (x_5,y_5) = (-96^3,\ &3^{12} 13^2 19^2 29\cdot31^2 37^2 \cdot53\cdot113\cdot173\cdot179)
\end{align*}
Since the absolute value obscures the polynomial's outputs, we fall back on an exhaustive search through the 16 different possibilities of choices for the signs of the $y_i$. Then using standard curve fitting we find that only one choice of signs, namely
$$(x_1,y_1),\ (x_2,y_2),\ (x_3,-y_3),\ (x_4,-y_4),\ (x_5,-y_5),$$
yields a monic polynomial. It is
\begin{align*}
M_{j(\tau)}(x) &= x^4+331531596 x^3-429878960946 x^2+109873509788637459 x+20919104368024767633\\
& = x^4 + 2^2 3^3 11\cdot29\cdot9623 x^3-2\cdot3^6 41\cdot1303\cdot5519x^2+3^{12} 103\cdot2007246533x+3^{15}17^323^329^3
\end{align*}
Since this is only a degree 4 polynomial over the reals, it is solvable by radicals and yields 4 roots. By comparing these to decimal approximations for $j(\tau)$ computed analytically we find  that 
$$j\left(\frac{1}{2}(1+i\sqrt{39})\right) = -\frac{27}{2} \left(6139474+1702799 \sqrt{13}+147 \sqrt{39 \left(89453213+24809858 \sqrt{13}\right)}\right).$$
Note that as a corollary we find that the algebraic norm $|j(\tau)| = 3^{15}17^323^329^3$. We could not have computed this from (\ref{GZ}) alone since the discriminants $-39$ and $-3$ are not relatively prime.
\end{subsection}

\end{section}


\begin{section}{Minimal polynomials of singular moduli from $\Xs_6$}\label{Chapter 4}
For the Shimura curve $\Xs_6$, the idea is essentially the same, but some modifications must be made. First, since there are no cusps on $\Xs_6$, there is no $\q$-expansion for the function $t: \Xs_6 \to \Pro^1$. This severely limits any analytical approaches. 

Second, the analog to (\ref{GZ}) in the Shimura curve case is considerably more complicated to compute. In \cite{JS}, Schofer attains an explicit formula for the average of a Borcherds form over CM points associated to a quadratic form of signature
$(n, 2)$. In the second half of \cite{JSt} Schofer shows that this generalizes the Gross-Zagier formula in the classical modular curve case. In \cite{Errthum} these techniques are applied to the Shimura curves $\Xs_6$ and $\Xs_{10}$ to algebraically compute the norms of CM points on these curves. We now further extend these methods in a manner similar to the previous section to compute the minimal polynomials of the singular moduli on $\Xs_6$.

Recall that any parametrization $t$ is invariant under a $\GL(\Z)$ action. Thus it suffices to assign its value at three CM points. Let $s_d$ denote the CM point of discriminant $d$. In \cite{Elk}, Elkies shows that the triangle group $\Gs$ is generated by $s_3$, $s_4$ and $s_{24}$ and makes the arbitrary, albeit natural, choice of defining the function's zeros and poles at those points. The computations in \cite{Errthum} follow suit, defining a map such that
\begin{equation}\begin{array}{c}
t(s_3)=\infty,\\
t(s_4)=0,\\ 
t(s_{24})=1,
\end{array} \label{Choices}\end{equation}
and determine $|t(s_d)|$, the rational norm of the singular moduli of $\Xs_6$.

We can now construct the minimal polynomial for $t(s_d)$ much as we did for $j(\tau)$ in the previous section.  In general, the algebraic degree of a singular modulus on a Shimura curve, $h_d$, can be computed using genus theory. For $\Xs_6$, there are exactly 27 rational singular moduli. Let $\zeta_d$ denote a rational singular moduli (e.g. $\zeta_d$ only makes sense for the 27 specific values of $d$ for which the singular moduli is, in fact, rational). Then we can define a new parametrization $t_d$ by making different choices than those in (\ref{Choices}). Specifically, if we choose instead, for $d \neq 3,4$,
\begin{equation}\begin{array}{c}
t_d(s_3)=\infty,\\
t_d(s_d)=0,\\ 
t_d(s_4)=\zeta_d,
\end{array} \label{AltCh}\end{equation}
it yields the relationship $$t_d(x) = \zeta_d - t(x).$$ Note, we can then formally define $t_4 = -t$.

Again, for a given $s_{d'}$ we can compute $|t_d(s_{d'})|$ in two different ways -- via the calculations in \cite{Errthum} and by definition:
\begin{align*}
|t_d(s_{d'})| &= \left|\prod_{\sigma} \sigma(t_d(s_{d'}))\right|\\
&=\left|\prod_{\sigma} \sigma(\zeta_d - t(s_{d'}))\right|\\
&=\left|\prod_{\sigma}  (\zeta_d -\sigma(t(s_{d'})))\right|&\text{(since $\zeta_d$ is rational)}\\
&=|M_{t(s_{d'})}(\zeta_d)|
\end{align*}

Repeating this for $h_{d'}+1$ choices of $r_d$ gives us sufficiently many points to exhaustively search the $2^{h_{d'}}$ possibilities for the monic polynomial $M_{t(s_{d'})}(x)$. Notice, though that since the calculations in \cite{Errthum} are not constrained by a property analogous to the relative primeness of discriminants, this method of calculation works for any singular modulus with $h_{d'} \le 26$, a much larger collection of points than in the classical case.

\begin{subsection}{Example: $t_6(s_{244})$}
We now consider the case of the CM point $s_{244} \in \Xs_6$. Genus theory shows that the image of $s_{244}$ is an algebraic number of degree 3. Since the norm calculator in the Shimura curve case does not have the relatively prime condition, we can use the algorithms in \cite{Errthum} to compute the following
\begin{eqnarray*}
|t_{4}(s_{244})| &=& \frac{2^6 3^{21}19^4}{17^6 29^6}\\
|t_{24}(s_{244})| &=& \frac{19^4 37^2 47^2 61}{17^2 29^6}, \\
|t_{40}(s_{244})| &=& \frac{3^{22} 83 \cdot 101\cdot 107 \cdot 163}{5^9 17^2 29^4}\\
|t_{52}(s_{244})| &=& \frac{2^{18} 3^{23} 37^2 103\cdot 131 \cdot 179 \cdot 199 \cdot 263}{5^{18}17^6 29^6}
\end{eqnarray*}
This gives the following 4 points on the curve $y=|M_{t(s_{d'})}(x)|$:
\begin{align*}
(x_1,y_1) = \Big(0,\ &\frac{2^6 3^{21}19^4}{17^6 29^6}\Big),\\
 (x_2,y_2) = \Big(1,\ & \frac{19^4 37^2 47^2 61}{17^2 29^6}\Big),\\
 (x_3,y_3) = \Big(-\frac{3^7}{5^3},\ & \frac{3^{22} 83 \cdot 101\cdot 107 \cdot 163}{5^9 17^2 29^4}\Big),\\
 (x_4,y_4) = \Big(\frac{2^23^7}{5^6},\ & \frac{2^{18} 3^{23} 37^2 103\cdot 131 \cdot 179 \cdot 199 \cdot 263}{5^{18}17^6 29^6}\Big),\\
\end{align*}
Using standard curve fitting, there are 8 different possibilities depending on the choices of sign for $y_i$. Only one choice of signs, namely
$$(x_1,-y_1),\ (x_2,y_2),\ (x_3,-y_3),\ (x_4,y_4),$$
yields a monic polynomial, which is
\begin{align*}
M_{t(s_{244})}(x) & = x^3-\frac{159511016412629892}{14357588953446649} x^2+\frac{2240284633411688496}{14357588953446649}x-\frac{87245036145162432}{14357588953446649}\\
& = x^3 - \frac{2^2 3^7 31\cdot 67\cdot 37223\cdot 235849}{17^6 29^6}x^2 + \frac{2^4 3^{14} 151\cdot 1187\cdot 163327}{17^6 29^6}x -  \frac{2^6 3^{21}19^4}{17^6 29^6}
\end{align*}
Although this is only a degree 3 polynomial over the reals and is thus solvable by radicals, there is no analytical approximation to compare to in order to determine exactly which root is $t(s_{244})$.
\end{subsection}

\end{section}


\begin{section}{Algebraic $abc$-ratios of singular moduli}\label{Chapter 5}

Oesterle and Masser's well-known $abc$-conjecture states that for relatively prime positive integers that satisfy $a+b=c$ there is a bound on how large $c$ can be in terms of the product of all the primes involved. More explicitly the conjecture asserts that given any $\epsilon>0$ there is a constant $C_\epsilon$ such that $$c \le C_\epsilon (\text{rad}(abc))^{1+\epsilon}$$ where $\text{rad}(n)$ is the product of all prime divisors of $n$.

Taking the constant as 1, one can measure the ``quality'' of an $abc$-example by the necessary size of $\epsilon$. For this reason we consider the $abc$-ratio
$$\alpha(a,b,c) = \frac{\ln(c)}{\ln(\text{rad}(abc))}.$$
To this date, the largest known $abc$-ratio is $$\alpha(2,3^{10}109,23^5) \approx 1.62991.$$
For comparison, the median value of $\alpha(a,b,c)$ for $1 \le a,b \le 100$ is approximately $0.429$ with the maximum being $\alpha(1,80,81) \approx 1.29203$. The standard benchmark has become that any example with $\alpha(a,b,c) > 1.4$ is considered ``good'' \cite{ABCweb}.

Vojta \cite{Vojta} generalized the $abc$-conjecture to number fields in the following way (as described in \cite{ABCweb}). Let $K$ be an algebraic number field and let $V_K$ denote the set of primes on $K$. Then any $v\in V_K$ gives an equivalence class of non-trivial norms on $K$ (finite or infinite). Let $||x||_v=|P|^{-v_P(x)}$ if $v$ is a prime defined by a prime ideal $P$ of the ring of integers $O_K$ in $K$ and $v_P$ is the corresponding valuation, where $|\cdot|$ is the absolute norm. Let $||x||_v=|g(x)|^e$ for all non-conjugate embeddings $g: K \to \mathbb{C}$ with $e=1$ if $g$ is real and $e=2$ if $g$ is complex. Define the height of any triple $a$, $b$, $c\in K^\times$ to be
$$H_K(a,b,c) = \prod_{v\in V_K} \text{max}(||a||_v, ||b||_v, ||c||_v),$$
and the radical (or conductor) of $(a,b,c)$ by
$$\text{rad}_K(a,b,c) = \prod_{P\in I_K(a,b,c)}|P|,$$
where $I_K(a,b,c)$ is the set of all prime ideals $P$ of $O_K$ for which $||a||_v$, $||b||_v$ , $||c||_v$ are not equal. Let $D_{K/\mathbb{Q}}$ denote the discriminant of $K$.

The (uniform) algebraic $abc$-conjecture then states that for any $\epsilon> 0$, there exists a positive constant $C_\epsilon$ such that for all $a$, $b$, $c \in K^\times$ satisfying $a+b+c=0$, we have
$$H_K(a,b,c) < C_{\epsilon}^{[K:\mathbb{Q}]}(|D_{K/\mathbb{Q}}|\text{rad}_K(a,b,c))^{1+\epsilon}.$$
Again, constraining the constant leads to a measure of the ``quality'' of an algebraic $abc$-example given by the algebraic $abc$-ratio
$$\gamma(a,b,c)=\frac{\ln(H_K(a,b,c)}{\ln(D_{K/\mathbb{Q}})+\ln(\text{rad}_K(a,b,c))}.$$
It's important to note that all of the above quantities are Galois-invariant, and thus an algebraic number and any of its conjugates have the same algebraic $abc$-ratio. Hence, it is sufficient to know only the minimal polynomial for the algebraic number.

The largest known algebraic $abc$-ratio is $$\gamma(w,(w+1)^{10}(w-1),2^9(w+1)^5) \approx 2.029229$$ where $w^2-w-3=0$. Since algebraic $abc$-ratios are in general slightly larger than $abc$-ratios, any example with $\gamma(a,b,c)>1.5$ is considered ``good'' \cite{ABCweb}.

\begin{subsection}{$abc$-ratios of Singular Moduli}
With the minimal polynomial available from the computations of Sections 3 and 4, one can now perform the necessary computations to find the algebraic $abc$-ratios of the singular moduli arising from the classical modular curve and the Shimura curve $\Xs_6$.

Figure 1 
\begin{figure}[t]
\centering\includegraphics[scale=.55]{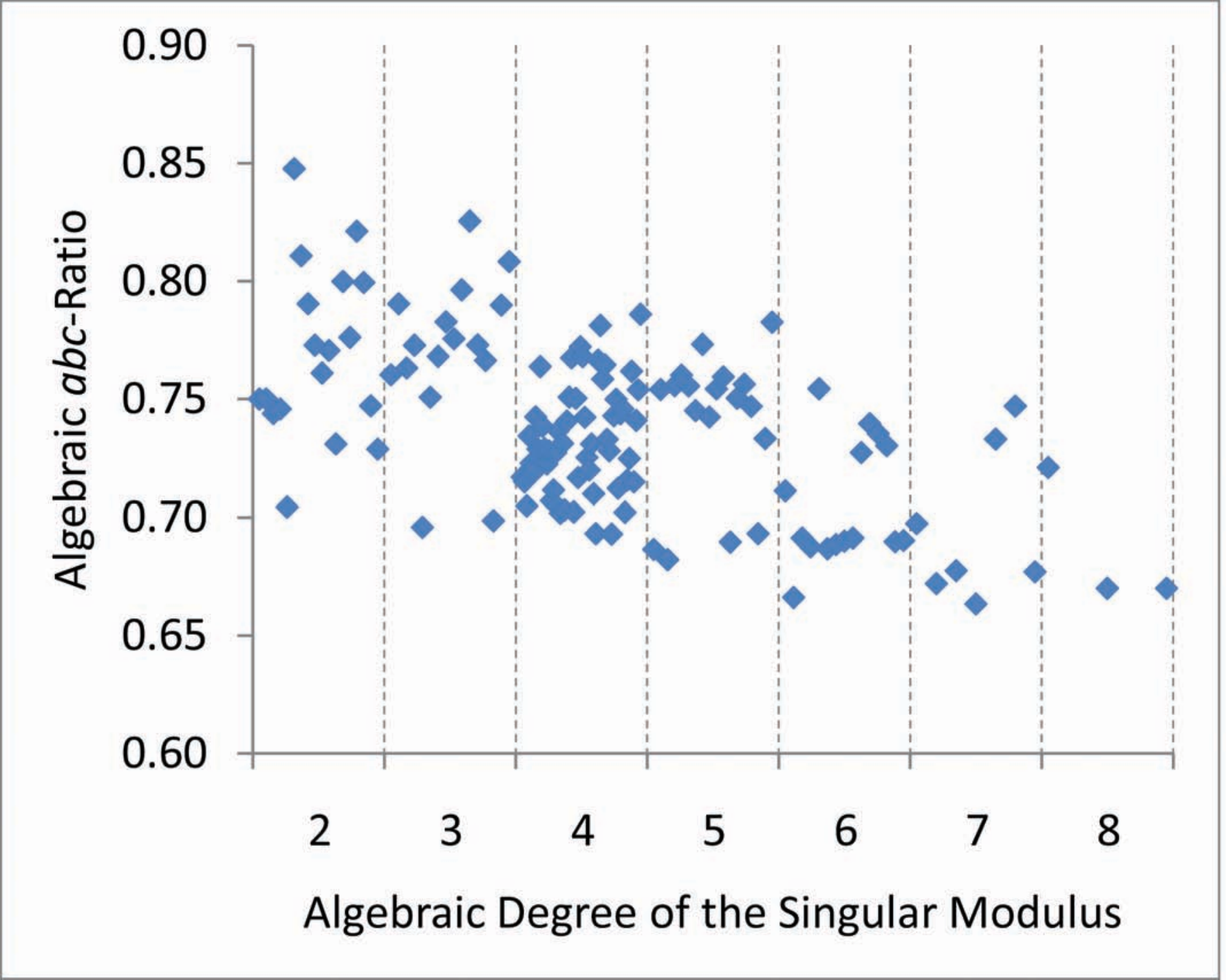}\caption{Algebraic $abc$-Ratio for Classical Singular Moduli}
\end{figure}
presents a plot of $\gamma(a,b,c)$ versus $[K:\mathbb{Q}]$ for the classical singular moduli. Figure 2 
\begin{figure}[t]
\centering\includegraphics[scale=.75]{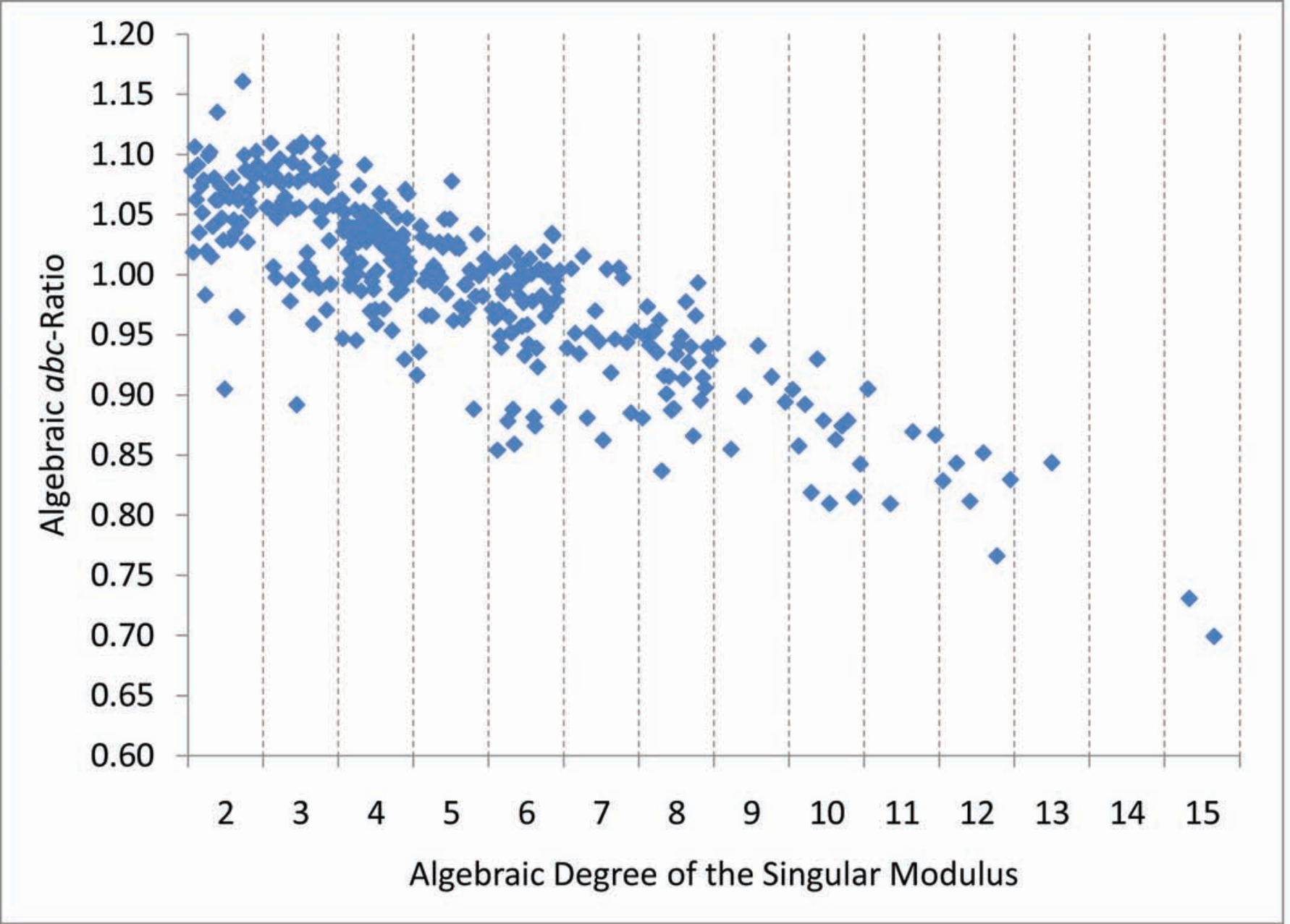}\caption{Algebraic $abc$-Ratio for $\Xs_6$ Singular Moduli}
\end{figure}
does the same for the Shimura curve $\Xs_6$. (Note: Data points with common degree were lexicographically ordered according to the coefficients of the minimal polynomial. Not all data points were available for the classical singular moduli due to software reaching computational limits. For $\Xs_6$, all data points of singular moduli up to the discriminant of 4744 are plotted.)

Unfortunately, in neither case are any ``good'' $abc$-examples found, even though the values are much higher than what would be typical of an arbitrary algebraic number. Lastly, the $abc$-ratios appear to be following a trend as the degree increases. It may be possible to use the general formulas in \cite{JS} to construct bounds on the algebraic $abc$-ratios of singular moduli.

\end{subsection}

\end{section}


\noindent {\bf Author Contact Information:}\\
Eric Errthum\\
Department of Mathematics \& Statistics\\
Winona State University\\
Winona, MN, USA 55987\\
eerrthum@winona.edu
\end{document}